
\magnification=1200
\baselineskip=14pt

\def\proof{\medskip\noindent{\bf Proof. }}

\def\O{{\cal O}}

\def\L{{\cal L}}

\def\P{{\bf P}}
\def\Q{{\bf Q}}
\def\P{{\bf P}}
\def\Z{{\bf Z}}
\def\C{{\bf C}}

\def\R{{\bf R}}

\def\f{\varphi}
\def\ra{\rightarrow}
\def\raa{\longrightarrow}
\def\iso{\simeq}

\def\har#1{\smash{\mathop{\hbox to .8 cm{\rightarrowfill}}
\limits^{\scriptstyle#1}_{}}}

\font\ninerm=cmr10 at 10truept
\font\ninebf=cmbx10 at 10truept
\font\nineit=cmti10 at 10truept
\font\bigrm=cmr12 at 12 pt
\font\bigbf=cmbx12 at 12 pt
\font\bigit=cmti12 at 12 pt

\def\smalltype{\let\rm=\ninerm \let\bf=\ninebf
\let\it=\nineit  \baselineskip=12pt minus .75pt \rm}

\def\bigtype{\let\rm=\bigrm \let\bf=\bigbf\let\it=\bigit
\baselineskip=14 pt minus 1 pt \rm}


\noindent February 1998

\vskip 2 cm

{\bigtype
\centerline{COHOMOLOGICAL INVARIANTS OF COMPLEX MANIFOLDS}
\medskip
\centerline{COMING FROM EXTREMAL RAYS}
\medskip
\centerline{Jaros{\l}aw A. Wi\'sniewski}
\medskip
\centerline{\it to Ma{\l}gosia, Jagna and Mary{\~n}a}
}


\bigskip

\beginsection Introduction: good rays and bad rays.

Topology has always been a very convenient tool for studying
compact complex manifolds. In particular, the cohomology ring
$H^*(X,\Z)$ of a compact complex manifold $X$ is a very important
invariant which provides a lot of information on the geometry of $X$.
In the present paper we use $H^*(X,\Z)$ to understand {\sl extremal ray
contractions} in the sense of the Minimal Model Program or Mori Theory.

Originally, an extremal ray contraction was invented to serve as a
step in the Minimal Model Program on the way to reach a minimal or
canonical model of a projective variety. However, it was soon observed
that such contractions can measure certain properties of {\sl
non-general} varieties. Recently the notion of an extremal ray and
its contraction has made its way to other branches of geometry:
non-projective geometry, Moishezon manifolds, deformation of complex
structure and symplectic geometry. At this point I should note
that in this paper I use the name ``extremal ray'' (and respectively
``extremal contraction'') for both {\sl Mori} (or {\sl Fano-Mori}) and
{\sl crepant} extremal rays, see the definitions in Section 2.

It is apparent from some of these applications of extremal rays that
not all rays can be equally treated. For example, in symplectic
geometry [19] (c.f. [17]) we distinguish the rays which
contain classes of {\sl good rational curves} with good
deformation properties.

In the present paper I try to explain this apparent non-equality of
extremal contractions on the level of cohomology. Namely, using the
invariants of $H^*(X,\R)$ which come from the Hard Lefschetz Theorem
I divide extremal rays in two classes: {\sl L-supported} and {\sl
L-negligible} (where ``L'' stands for ``Lefschetz''). Roughly
speaking: L-supported rays are strongly distinguishable in topology
while L-negligible rays have very mild geometry. Now the decision
which rays are good and which are bad depends on the taste, or
rather on a possible purpose they should serve for.

Let me summarize the properties of these two classes of rays.  Each
L-supported ray $R$ defines a hyperplane $R^\perp\subset H^2(X,\R)$
on which Lefschetz duality degenerates (i.e.~$R^\perp$ is a component
of the {\sl Lefschetz discriminant} defined in Section 1), so the
hyperplane is a cohomology ring invariant. The hyperplane $R^\perp$
carries a multiplicity $m(R)$ (a cohomology ring invariant) which is
related to the geometry of the ray $R$, for example $m(R)\geq l(R)-1$
where $l(R)$ is the length of the ray $R$ defined in [21].  The
number of L-supported rays is finite and it can be bounded (counting
them with their multiplicities) by the number $\sum_{k=1}^{n}k\cdot
b_{n-k}(X)$.

Although the number of L-negligible rays may be infinite and they are
invisible in the cohomology ring, their geometry is easier than that
of L-supported rays. For example, in dimension 2 and 3 the only Mori
extremal rays come from blowing up a smooth variety along a
codimension 2 smooth center.  Further, in codimension 4 contractions
of Mori extremal rays can be well described (Theorem 3.4).  Moreover,
each L-negligible ray contains lots of ``good'' rational curves whose
deformation is of the expected dimension. In effect, L-negligible
rays are invariant under deformations of complex structure and can be used
to compute Gromov-Witten invariants in symplectic geometry.

The paper is organised as follows. After reviewing basic
definitions concerning topology and Mori Theory, in Section 1 some
consequences of Hard Lefschetz Theorem are worked out. In particular,
the notions of {\sl Lefschetz discriminant} and {\sl L-supported
homology class} are introduced. In Section 2, after recalling some
further facts from Mori Theory concerning extremal contractions a
{\sl representability theorem} is proved; it allows to represent
cohomology in the hyperplane $R^\perp\subset H^2(X,\R)$ by cohomology
of the target of the contraction of $R$.  The geometry
of L-negligible rays is studied in Section 3.

The present paper was motivated by some recent publications on
extremal rays in different branches of geometry. Apart of [19] and
[17] I should mention also [20], [18] and [2]. I would like to thank
the authors of the last two papers for their kind sending me their
preprints.  In some sense the present paper is a sequel of [23] where
the idea of using Hard Lefschetz Theorem in the context of Mori
Theory was born.

During the prepation of this paper I had a fellowship at the
Institute of Mathematics of Polish Academy of Sciences.  In April and
May 1997 I lectured on the topology of complex manifolds at The
University of Trento where I presented the results of this paper.  I
am very grateful to the participants of these lectures, especially to
Marco Andreatta, Edoardo Ballico and Massimiliano Mella, for their
remarks and criticism which corrected some of my wrongthinking.  I
was also partially supported by Polish KBN.  I would like to thank
all the above institutions for their kind support.


\beginsection Notation and definitions.

Throughout the present paper $X$ is a complex manifold of complex
dimension $n$ and moreover --- with the exception of a local set-up in
Section 2 --- $X$ is compact or just projective. We consider the
complex topology of $X$.  We will usually work with homology and
cohomology of $X$ with coefficients in the field of real numbers
$\R$; a cycle (or cocycle) with real coefficients will be called
integral or rational if it is coming from homology (or cohomology)
with integral or, respectively, rational coefficients.

A great deal of information about the relation of the complex structure
of $X$ with its topology is carried by the exponential sequence on $X$:
$$0\raa\Z_X\raa\O_X\har{exp(2\pi i\cdot\phantom{x})}\O^*_X\raa 0$$
the long cohomology sequence of which defines the link betweeen the
Picard group of $X$ and its second cohomology group and in particular
the first Chern class map:
$$ H^1(X,\Z_X)\raa H^1(X,\O_X)\raa PicX = H^1(X,\O^*_X)
\har{c_1} H^2(X,\Z_X)\raa H^2(X,\O_X)$$
The image of $c_1$ in the torsion-free cohomology $H^2(X,\Z)/({\rm
torsion})$ coincides with the group of divisors modulo
numerical equivalence $PicX/\equiv$.

Let $N^1(X)\subset H^2(X,\R)$ and $N_1(X)\subset H_1(X,\R)$ be $\R$-linear
subspaces spanned by, respectively, cohomology and homology classes of,
respectively, holomorphic divisors and curves on $X$. In other words
$N^1(X)$ is spanned by the image of the map $c_1: PicX\ra H^2(X,\Z)$
composed with the extension of coefficients $H^2(X,\Z)\ra H^2(X,\R)$.

The topological intersection of cycles and cocycles restricts to $N_1(X)$
and $N^1(X)$ and coincides with the intersection product which can be
defined in algebro-geometric set up; in both cases the intersection of a
cocycle $\chi$ with a cycle $\alpha$ will be denoted by $\chi\cdot\alpha$.
The intersection product gives a nondegenerate pairing on $H^2(X,\R)\times
H_2(X,\R)$ and $N^1(X)\times N_1(X)$ and thus we will frequently identify
any space in question with the dual of its pairing partner. For a given
non-zero cycle $\alpha\in H_2(X,\R)$ we define the perpendicular
hyperplane $\alpha^\perp:=\{\chi\in H^2(X,\R): \chi\cdot\alpha=0\}$;
similarly we define the hyperplane perpendicular to a non-zero cocycle.

Inside the above spaces we consider the following cones. The cone of
curves ${\cal C}\subset N_1(X)$ and the cone of nef divisors ${\cal
P}\subset N^1(X)$ are $\R^*_{>0}$-spanned on, respectively, the classes of
curves and numerically effective divisors, i.e. ${\cal P}:= \{\chi\in
N^1(X): \forall \alpha\in{\cal C} \ \chi\cdot\alpha\geq 0\}$.  That is
${\cal P}={\cal C}^\vee$ in the sense of the intersection pairing of
$N^1(X)$ and $N_1(X)$. Dually, $\bar{\cal C}={\cal P}^\vee$, where
$\bar{\cal C}$ denotes the closure of ${\cal C}$. We recall that, by the
Kleiman criterion of ampleness, a line bundle $\L$ is ample if $c_1(\L)$
is in the interior (in the sense of the topology on $N^1(X)$) of the cone
${\cal P}$, which we will denote by ${\cal P}'$.


\beginsection 1. Hard Lefschetz Theorem.

In the present section we deal with the classical result of
Lefschetz.  Our purpose is to understand its impact on the structure
of the cohomology ring of a complex projetive manifold and on the
second cohomology in particular. For the exposition of the result the
reader may consult the classical textbook [8] or a recent excellent
survey [Looijenga].

Let $X$ be a compact complex manifold of dimension $n$. In the present
section we discuss some properties of the cohomology ring $\bigoplus_m
H^m(X,\R)$ which is equipped with the cup product $\cup$.
We recall that $b_m(X)=dim_\R H^m(X,\R)$ is the $m$-th Betti number of
$X$.

\proclaim Definition. We will say that $\eta\in H^2(X,\R)$ satisfies
Lefschetz condition if the $k$-th cup product map $$
\matrix{L_k(\eta):& H^{n-k}(X,\R)&\raa& H^{n+k}(X,\R)\cr
&\chi&\raa&\chi\cup\eta^{\cup k}} \leqno (\dag)$$ is an isomorphism
for $k=1,\dots n$ (then we will call it {\sl Lefschetz duality}).

By abuse, if $\eta=\eta(D)$ is the cohomology class of
a Cartier divisor $D$ or, equivalently, the first Chern class $c_1(\L)$ of
the associated line bundle $\L={\cal O}_X(D)$ then we will say that $D$ and
$\L$ satisfy the Lefschetz condition as well. The above definition is
motivated by

\proclaim Hard Lefschetz Theorem. If $\L$ is an ample line bundle then it
satisfies the Lefschetz condition.

Alternatively, instead of considering the map $L_k(\eta)$ we can take
$(-1)^{n-k}$-symmetric bilinear pairing
$$\matrix{H^{n-k}(X,\R)\times H^{n-k}(X,\R)&\raa& \R\hfill\cr
(\chi,\nu)&\raa&
\chi\cup\nu\cup\eta^{\cup k}}
\leqno(\dag\dag)$$
In other words, the choice of $\eta$ gives a 2-form $A_k(\eta)$ which
lives in $S^2(H^{n-k}(X,\R)^*)$ or $\bigwedge^2(H^{n-k}(X,\R)^*)$ ---
depending on whether $n-k$ is even or odd --- and which is
nondegenerate if and only if $\eta$ satisfies the $k$-th cup product
map is an isomorphism. In a fixed basis the form $A_k(\eta)$ can be
represented as a symmetric (resp.~skew-symmetric) $b_{n-k}\times
b_{n-k}$ matrix, the coefficients of which are homogeneous of degree
$k$ in $\eta$. Therefore, since $A_k(\eta)$ is nondegenerate if and
only if $detA_k(\eta)\ne 0$, we obtain the following

\proclaim Lemma 1.1. Assume that $X$ is a projective manifold.
Let $\Delta_k=\Delta_k(X):= \{\eta\in H^2(X,\R): L_k(\eta)$ is not an
isomorphism$\}$. Then $\Delta_k$ is the zero set of a homogeneous
polynomial $\delta_k=det A_k$ of degree $kb_{n-k}$ (if $n-k$ is odd then
$\delta_k$ is the square of the pfaffian of $A_k$).  Moreover, identifying
rational homology $H_2(X,\Q)$ with rational linear forms on $H^2(X,\R)$ we
have $\delta_k\in S^{kb_{n-k}}(H_2(X,\Q))$.

\proof Lemma follows from the preceeding discussion in which the field
$\R$ can be replaced by $\Q$ so that $\delta_k\in
S^{kb_{n-k}}(H_2(X,\Q))$.
\medskip

The set $\Delta_k\subset H^2(X,\R)$ which is introduced just above will
be called the $k$-th Lefschetz discriminant. The immediate consequence of
Lefschetz theorem is that for $k=1,\dots n$ the Lefschetz discriminant
$\Delta(X)=\bigcup\Delta_k$ does not meet rational points in ${\cal P}'$.
\medskip

\noindent {\bf Example.}
Let $X=\C^n/\Gamma$ be a compact complex torus. The cohomology
ring of $X$ can be identified with the ring of alternating real forms on
$\C^n=\R^{2n}$ that is $$\bigoplus_k H^k(X,\R)\iso \bigoplus_k
\bigwedge^k(\R^{2n}).$$ It is easy to verify that a form $\eta$
satisfies Lefschetz condition if and only if $\eta^{\wedge n}\ne 0$,
i.e.~if it is symplectic. Indeed, this is equivalent that in some
basis it can written as $\eta= dx_1\wedge dx_2+dx_3\wedge dx_4+\dots
+dx_{2n-1}\wedge dx_{2n}$ and thus, for a suitable choice of the
complex structure on $\R^{2n}$ we have $\eta=\sqrt{-1}\cdot
(dz_1\wedge d\bar z_1+dz_2\wedge d\bar{z}_2+\dots +dz_n\wedge
d\bar{z}_n)$. In this case the Lefschetz condition is defined by the
$n$-th product map, or equivalently for $i=1,\dots n$ we have
$\Delta_i\subset
\Delta_n$.

\proclaim Definition. A non-zero homology class $\alpha\in H_2(X,\R)$ is
called Lefschetz supported (or simply L-supported) if it divides
$\delta_k$ for some $k=1,\dots n$.

By abuse, a line $\R\cdot\alpha\subset H_2(X,\R)$ or a ray $\R_{\geq
0}\cdot\alpha\subset H_2(X,\R)$ will be called L-supported if $\alpha$ is
L-supported. Moreover, we will say that $\alpha$ is of type $k$ and
multiplicity $m$ if it divides $\delta_k$ with multiplicity exactly $m$.
In this notation we have the following corollary to the previous Lemma.

\proclaim Corollary 1.2. Let $X$ be a projective manifold.
Then the  set of L-supported lines in $H^2(X,\R)$ is finite.
More exactly, if $\#(k,m)$ is the number of lines of type $k$ and
multiplicity $m$ then for $k=1,\dots n$ we have
$$\sum_{m\geq 1}m\cdot\#(k,m)\leq k\cdot b_{n-k}.$$

The reason to introduce the notion of L-supported homology class is the
following consequence of Lefschetz Theorem.

\proclaim Proposition 1.3. Assume that $X$ is a projective manifold.
Let $\alpha$ be an L-supported rational homology class in
$H_2(X,\R)$. Then
\item{(1)} the cone ${\cal P}$ is not cut by the hyperplane $\alpha^\perp$,
or equivalently, either $\alpha\cdot{\cal P}\geq 0$ or $\alpha\cdot{\cal
P}\leq 0$;
\item{(2)} if $\beta\in N_1(X)$ is the projection of $\alpha$ along
$N^1(X)^\perp$ (i.e. $\beta$ is the unique class in $N_1(X)$ such that
$(\beta-\alpha)\cdot N^1(X)=0$) then either $\beta\in\bar{\cal C}$ or
$-\beta\in\bar{\cal C}$.

\proof Immediate.

\medskip

\noindent {\bf Remark.}
Lefschetz condition is $\R^*$ invariant, that is $\eta\in H^2(X,\R)$
satisfies it if and only if $a\eta$ does, for any $a\in\R\setminus
\{0\}$. Although, for a given complex structure the notion of ample cone
is only $\R^*_{>0}$ invariant, we note however that the topological
invariants should not distinguish an ample divisor from its opposite.

Indeed, the topological space of a given complex manifold $X$
supports also its conjugate structure $\bar X$ which is associated to
the conjugation $\bar{\phantom{x}} :\C\ra\C$ of the base field.  If
$\eta\in H^2(X,\Z)$ is a class of an ample divisor on $X$ then
$-\eta$ is the class of an ample divisor on $\bar X$.


\beginsection 2. Extremal rays and extremal contractions.

The task of the present paper is to apply the formalism which we have
introduced in the previous section in the situation appearing in Mori
theory of extremal ray contractions. In the present section we use
the language and the fundamental results of the Minimal Model
Program. For an introduction and an exposition to the program we
suggest [5] or more advanced [12]. The cone theorem (smooth case) is
explained in [16]. The local analytic version of the contraction
theory which we will rely on is in [10].

A {\sl contraction} is a proper surjective map $\f:X\ra Y$ of normal
irreducible varieties with connected fibers such that
$\f_*\O_X=\O_Y$. We assume that a contraction is not an isomorphism.
The map $\f$ is birational or otherwise $dimY< dimX$, in the latter
case we say that $\f$ is of fiber type.  The exceptional locus
$E(\f)$ of a birational contraction $\f$ is equal to the smallest
subset of $X$ such that $\f$ is an isomorphism on $X\setminus E(\f)$.
The contraction $\f$ is called {\sl Fano-Mori} (or just {\sl Mori
contraction}) if the anti-canonical divisor $-K_X$ is $\f$-ample.  If
the map $\f$ is birational and $K_X=\f^*K_Y$ then we say that $\f$ is
{\sl crepant}.

The contraction $\f$ yields two maps $\f^*: H^2(Y,{\bf Z})\ra
H^2(X,{\bf Z})$ on cohomology and also $\f^*:N^1(Y)\ra N^1(X)$, which
we denote similarly.  We say that $\f$ is {\sl elementary} if
$dim(N^1(X)/\f^*N^1(Y))=1$.

If both $X$ and $Y$ are projective then $\f^*({\cal P}(Y))$ is a face
of the cone ${\cal P}(X)$. Moreover the contraction $\f$ kills
(i.e.~contracts to points) the {\sl holomorphic} curves whose classes
are perpendicular to the face $\f^*({\cal P}(Y))$. In particular, an
elementary contraction defines a 1-dimensional face (a ray) in ${\cal
C}(X)$. The ray is called {\sl Mori ray} or {\sl crepant ray} if its
contraction is of the respective type.

It is remarkable that in the situation covered by the Minimal Model
Program the passage from contractions to rays (or faces) of ${\cal
C}(X)$ can be reversed. That is, if a projective normal variety $X$
has suitable singularities (e.g.~it is smooth) then a ``good''
face of the cone ${\cal C}(X)$ admits a contraction.  In particular,
if $R\subset {\cal C}(X)$ is a ray such that for some $D\in
R^{\perp}\cap {\cal P}(X)$, $D^\perp\cap{\cal C}(X)=R$ and either
$R\cdot K_X<0$ or $R\cdot K_X=0$ and $D^n>0$, then there exists a
contraction $\f_R:X\ra Y_R$ of $R$, where $Y_R$ is a projective
normal variety. The divisor $D$ is then called a {\sl good supporting
divisor} for $R$ and it is a pullback of an ample divisor from $Y_R$.
The contraction $\f_R$ is either Fano-Mori or crepant depending on
whether $R\cdot K_X<0$ or $R\cdot K_X=0$, respectively.  This is just
a version of

\proclaim Kawamata-Shokurov Contraction (or Base-Point-Free) Theorem.
Assume that $X$ is a projective manifold. Let $D\in PicX$ be a
line bundle such that $[D]$ as well as $[aD-K_X]$ are in ${\cal
P}(X)$ for some $a>0$.  If moreover $(aD-K_X)^n>0$ then there exists
a contraction $\f:X\ra Y$ such that $D=\f^*(D_Y)$ for some ample
divisor $D_Y\in PicY$.

We note that both Fano-Mori and crepant case are covered by the same
theorem.  For this reason, in the present paper frequently we do not
make any distinction between Fano-Mori or crepant elementary
contractions --- we will call both {\sl extremal contractions} which
is a slight abuse of the usual notation.

\medskip

Our task is to prove that in this range the word ``elementary'' can
be used in a broader sense, i.e.~any extremal contraction is
elementary also topologically. That is $N^1(X)/\f^*(N^1(Y))\iso \R$
\ if\  $\f:X\ra Y$ is an extremal contraction.

First we need a local version of the Base-Point-Free Theorem.
Namely, in the following lemma we will use the set-up and the
language explained in [10] pp.~102--106. That is, $X$ is a complex
manifold and the contraction $\f:X\ra Y$ is now a projective morphism
onto a normal analytic space $Y$ with a fixed point $y\in Y$.
Shrinking $Y$ if necessary we may assume that it satisfies the
assumptions 1.11 of [ibid]. Now we adapt the definitions of
$\f$-nefness and bigness for the map $\f: (X,\f^{-1}(y))\ra (Y,y)$ as
explained at pp.~105--106 of [ibid].

\proclaim Lemma 2.1.
Let $\f:(X,\f^{-1}(y))\ra (Y,y)$ be as above.
Assume that $-K_X$ is $\f$-big and nef with respect to $\{y\}$.
If $L\in Pic X$ is a line bundle such that $L\cdot C=0$ for any curve
$C\subset \f^{-1}(y)$ then there exists an open neighbourhood
$Y'$ of $y$ such that $L$ is trivial on $X'=\f^{-1}(Y')$,
i.e.~$L_{|X'}\iso\O_{X'}$.

\proof Because of Theorem 1.3' in [10]
$L^{\otimes m}$ is $\f$-spanned on $\f^{-1}(y)$ for $m\gg 0$.
Thus, since $L$ is numerically trivial on $\f^{-1}(y)$,
we can choose a section of $L^{\otimes m}$ which does not
vanish anywhere on $\f^{-1}(y)$. Therefore, $L^{\otimes m}$ is
trivial in a neighbourhood of $\f^{-1}(y)$ for $m\gg 0$,
so is $L$ itself.

\proclaim Corollary 2.2.
Let $\f:X\ra Y$ be a contraction of a projective manifold.
Suppose that $-K_X$ is $\f$-big and nef.
Then for any $y\in Y$ and any non-trivial $L\in (R^1\f_*\O_X^*)_y$
there exists a holomorphic curve $C\subset \f^{-1}(y)$ such that
$L\cdot C\ne 0$.  Thus the group $(R^1\f_*\O_X^*)_y$ is a finitely
generated abelian group with no torsion.

\proof The first part is just the previous Lemma. To get the second 
part we use the intersection product $(R^1\f_*\O_X^*)_y \times
H_2(\f^{-1}(y),\Z)\ra \Z$ and the first part to get
$R^1\f_*\O_X^*\hookrightarrow Hom(H_2(\f^{-1}(y),\Z),\Z)$.

\medskip

Now we are ready to prove the main technical result of this part.

\proclaim Proposition 2.3.
Let $\f:X\ra Y$ be a contraction of a projective manifold.  Assume
that $-K_X$ is $\f$-big and nef.  Then $R^1\f_*\Z_X=0$, moreover the
cohomology map $\f^*: H^2(Y,\Z)\ra H^2(X,\Z)$ is injective and the
cokernel of the natural inclusion $Pic(X)/\f^*Pic(Y)\hookrightarrow
H^2(X,\Z)/\f^*H^2(Y,\Z)$ is a torsion group.

\proof Let us consider the direct image of the exponential sequence on
$X$: $$0\raa R^1\f_*\Z_X\raa R^1\f_*\O_X\raa R^1\f_*\O_X^*\raa
R^2\f_*Z_X\raa R^2\f_*\O_X$$ Because of the vanishing of
$R^i\f_*\O_X$ for $i\geq 1$ (see e.g.~[10], Theorems 1.2 and
1.2') we get $R^1\f_*\Z_X=0$ and isomorphism of sheaves of abelian
groups $R^1\f_*\O^*_X\iso R^2\f_*\Z_X$.  We consider a commutative
diagram with exact rows coming from Leray spectral sequence and
vertical arrows coming from cohomology of the exponential sequence
\par
$$\matrix{
\hfill 0&\raa&PicY &\har{\f^*}&PicX&\har{u}&H^0(Y,R^1\f_*\O^*_X)\cr
& &\downarrow& &\downarrow&&\downarrow\cr
H^0(Y,R^1\f_*\Z_X)=0&\raa&H^2(Y,\Z)&\har{\f^*}&H^2(X,\Z)&\har{v}
&H^0(Y,R^2\f_*\Z_X)
}$$
\par\noindent
Because of 2.2 the most-to-the-right vertical arrow in an isomorphism
of torsionfree abelian groups which will use to identify them.
Let $U:=im(u)$ and $V:=im(v)$. Then, both $U$ and $V$ are finitely generated
abelian with no torsion and because of this identification
we can write $U\hookrightarrow V$.

Let us consider a free abelian group generated by curves contracted
by $\f$, that is ${\cal Z}_1(\f):=\{\sum a_iC_i: a_i\in \Z\}$ where
$C_i\subset X$ are homolorphic curves contracted by $\f$, only finite
number of $a_i$ is non-zero and ${\cal Z}_1(\f)$ has a natural group
structure. We have a natural intersection product
$$H^0(Y,R^1\f_*\O^*_X)\times {\cal Z}_1(\f)\raa\Z.$$ In view of 2.2,
for any non-zero $L\in H^0(Y, R^1\f_*\O_X^*)$ there exists $Z\in
{\cal Z}_1(\f)$ such that $L\cdot Z\ne 0$.
Therefore the product
yields a map ${\cal Z}_1(\f)\ra Hom(V,\Z)$ whose cokernel is a torsion
group.

Note that to any 1-cycle $Z\in{\cal Z}_1(\f)$ we can associate its
class $[Z]_X\in H_2(X,\Z)$.  Now $Z\cdot U=0$ if and only if $[Z]_X$
is zero on $c_1(PicX)\subset H^2(X,\Z)$. However, since $Z$ is a
1-cycle of holomorphic curves it follows that actually $[Z]_X$
vanishes on the whole cohomology $H^2(X,\Z)$ hence $Z\cdot V=0$.
Thus the image of ${\cal Z}_1\ra Hom(V,\Z)$ does not meet non-zero
elements in
$ker(Hom(V,\Z)\ra Hom(U,\Z))\iso Hom(coker(U\ra V),\Z)$
so the latter group is zero which implies that $coker(U\hookrightarrow V)$
is a torsion group.
\medskip

\noindent{\bf Example.} The quotient $V/U$ may be actually
non-zero (notation as in the proof). Indeed, let $Y$ be an abelian
surface and $\f: X\ra Y$ a $\P^1$-bundle in complex topology which
does {\sl not} come from a rank 2 vector bundle on $Y$. Then
$R^2\f_*\Z_X\iso R^1\f_*\O^*_X\iso\Z_Y$, $H^3(Y,\Z)\iso\Z^4$ and
using further terms of the Leray spectral sequence we get a diagram
$$\matrix{
H^1(X,\O^*_X)&\har{u}&H^0(Y,R^1\f_*\O^*_X)&\iso&\Z&\raa&H^2(Y,\O^*_Y)&
&\cr
\downarrow   &       &\downarrow          &    &  &    &\downarrow   &   &\cr
H^2(X,\Z)    &\har{v}&H^0(Y,R^2\f_*\Z_X)  &\iso&\Z&\raa&H^3(Y,\Z)    &\iso&\Z^4
}$$
Thus the map $v$ is surjective while $u$ is onto $2\Z\subset\Z$ since
$\O(1)$ on the fiber does not extend to $X$.
Finally, let us note that from the above sequence it follows that 
$V/U\subset H^2(Y,\O^*_Y)$ even if $\f$ is not a projective bundle.
In fact the group $V/U$ seems to be a good generalization of the
invariant $\delta_r$ which is defined for projective bundles as in 
[6].
\medskip
As a corollary we get the following {\sl representability theorem}

\proclaim Theorem 2.4.
Let $\f_R: X\ra Y$ be an extremal contraction of a ray $R$
of a  projective manifold. Then $R^\perp=\f_R^*(H^2(Y,\R))$.

In other words any {\sl topological} cocycle perpendicular to $R$ is
represented by a pullback of a cocycle from the target of the
contraction $\f_R$.


\beginsection 3. L-negligible extremal rays.

In this section we assume that $\f=\f_R:X\ra Y$ is an extremal
contraction of a ray $R$ on a smooth projective variety of dimension
$n$. If $R$ is L-supported then it is noticeable from the topological
view point as it gives a strong trace in the cohomology by giving a
component of the Lefschetz discriminant. Moreover, we have seen that
the number of L-supported rays on $X$ is finite. On the other hand,
although from the Mori Cone Theorem in [16] it follows that the
extremal rays in the cone ${\cal C}(X)$ are dicrete, its total number
on $X$ may be infinite.  (A simple example is a $\P^2$ blown up in 9
points so that the resulting surfaces has infinite number of
$(-1)$-curves, see e.g.~[5] example 4.6.4)

Therefore one may be tempted to conclude that ``{\sl a generic}
extremal ray on a variety $X$ is not L-supported''.  This gives a
motivation for understanding extremal rays which are not
L-supported; we will call them {\sl L-negligible}, for short.

\proclaim Proposition 3.1.
Let $\f_R: X\ra Y$ be an extremal contraction of a
ray $R$ of a projective manifold.  If there exists a subset
$S\subset X$ such that $2dim_\C S-dim_\C\f(S) \geq n+m$ for
some positive integer $m$ then $R$ is L-supported of type
$k:=2dim_\C S-n$ and multiplicity $\geq m$.

\proof Let us set $dim_\C S=a$, $dim_\C \f(S)=b-1$ so that $k=2a-n$
and $m=k-b+1$.
Let $[S]_X\in H_{2a}(X,\R)$ and $\nu_S\in H^{2n-2a}(X,\R)$ denote the
homology and cohomology class of $S$, respectively.
Consider a cocycle $\eta_Y\in H^2(Y,\R)$.
We claim that $\nu_S\cup\f^*(\eta_Y)^{\cup b}=0$ in $H^{2n-2a+2b}(X,\R)$.
This is equivalent to show that $[S]_X\cap\f^*(\eta_Y)^{\cup b}=0$ in
$H_{2b-2a}(X,\R)$. Let $i: S\ra X$, $j:\f(S)\ra Y$ be embeddings and
$\f_S: S\ra \f(S)$ the restriction of $\f$, so that $\f\circ
i=j\circ\f_S$. Then $$[S]_X\cap\f^*(\eta_Y)^{\cup b}=i_*([S]\cap
i^*(\f^*(\eta_Y))^{\cup b})= i_*([S]\cap\f^*_S(j^*(\eta_Y)^{\cup
b}))$$ but $j^*(\eta_Y)^{\cup b}\in H^{2b}(\f(S),\R)=0$, because
$dim_\C\f(S)=b-1$.

Let us choose a cocycle $\chi_0\not\in R^\perp$ and consider an
arbitrary $\chi\in H^2(X,\R)$ which we can write as
$\chi=\f^*(\eta_Y)+t\chi_0$ for an appropriate choice of $t\in\R$ and
$\eta_Y\in H^2(Y,\R)$.  Namely, we set
$t=(\alpha\cdot\chi)/(\alpha\cdot\chi_0)$ for a nonzero $\alpha\in R$
and then $\chi-t\chi_0\in R^\perp$ is in the image of $\f^*$ by the
representability.  Then $$\chi^{\cup k}\cup\nu_S= t^{k-b+1}{k\over
b-1}\cdot\eta_Y^{b-1}\cup\chi_0^{k-b+1}\cup\nu_S+
\dots + t^k\chi_0^k\cup\nu_S.$$
Therefore, if we write the matrix of the Lefschetz duality form
$A_k(\chi)$ in a basis $\{\nu_0:=\nu_S,\nu_1,\dots\}$ of the
cohomology $H^{n-k}(X,\R)$ then the first row is divisible by $t^m$.
Thus $\delta_k$ is divisible by $t^m$ which concludes the proof of
the proposition.

\bigskip

The above result should be compared with the following property of
extremal contractions, see [9] Theorem 0.4, [7] Lemma 2.5 and
[22] Theorem 1.1.

\proclaim Theorem 3.2.
Let $\f=\f_R: X\ra Y$ be an extremal contraction of a ray $R$ and let
$S\subset E(\f)$ be an irreducible component of the exceptional locus
of $\f$.  We define the length $l_S(R)$ of $R$ at $S$
$l_S(R):=min\{-K_X\cdot [C]\}$ where $C$ is a rational curve
passing through a general point of $S$ and $[C]\in R$.  Then
we have
$$2dim_{\C}S-dim_{\C}\f(S)\geq n+l_S(R)-1.$$

The inequality appearing in the above theorem is sometimes referred
to as {\sl fiber-locus inequality}. The bounds provided by the
fiber-locus inequality and the inequality from 3.2 give a narrow
space for L-negligible extremal contractions. If $R$ is crepant then
$l_S(R)=0$ and $2dim_{\C}S-dim_{\C}\f(S)$ is either $n-1$ or
$n$.

If however $R$ is a Mori ray then the exceptional locus of $\f_R$ is
covered by rational curves such that $-K_X\cdot [C]=1$. Moreover, for
any component $S$ of $E(\f)$ we have $2dim_{\C}S-dim_{\C}\f(S)=n$.
As the consequence, Mori L-negligible rays can be pretty
well described. Thus for the rest of this section we assume that $R$
is a Mori L-negligible ray.  A similar situation (when the fiber-locus
inequality becomes an equality) was considered in Lemma 1.1 of
[4].

\proclaim Lemma 3.3.
Let $\f:X\ra Y$ be a contraction of a Mori L-negligible ray.
For an irreducible component $S\subset E(\f)$
let  $r:=codim_X S=dim_{\C}S-dim_{\C}\f(S)$.
If $F$ is an irreducible component of a fiber of $\f_{|S}$ such that
$dim_{\C}F=r$ then its normalization $f:\tilde F\ra F\subset S$
is isomorphic to a projective space  $\tilde F\iso\P^r$ and
$f^*\O_X(-K_X)\iso\O(1)$.

The dimension of the exceptional locus and the dimension of fibers of
a Fano-Mori contraction of a L-negligible ray are tied up very closely so
that in small dimensions (or small codimensions) fibers of such
contractions are small and thus were thoroughly studied. The
subsequent structure theorem summarizes some of the known properties
of such contractions. We refer the reader to [3] for an exhausting
description of Fano-Mori contraction with fibers of dimension $\leq
2$.

\proclaim Theorem 3.4.
Let $\f=\f_R:X\ra Y$ be a Fano-Mori contraction
of an L-negligible ray on a smooth projective variety.
Let $E=E(\f)$ denote the exceptional locus of $\f$
and let $Z:=\f(E)\subset Y$ be the exceptional locus of $\f^{-1}$.
If $Z_k:=\{y\in Y: dim_{\C}f^{-1}(y)\geq k\}\subset Z$ then
$dim_{\C}Z_k\leq n-2k$. Moreover,
\item{(i)} if $dim_{\C}E=n-1$ then $E$ is an irreducible divisor,
$dim_{\C}Z=n-2$; outside the set $Z_{2}$
both $Y$ and $Z$ are smooth and
outside of the set $Z_3$ the map $\f^{-1}$ is a blow-up of $Y$ along $Z$.
\item{(ii)} if $dim_{\C}E\leq n-2$ then $Z=Z_2$ and outside a set
$Z'\subset Z$, $dim_{\C}Z'\leq n-5$, the map $f^{-1}$ is a small resolution
of a family of Veronese cone singularities. That is, for any
$y\in Z\setminus Z'$ there exists a neighbourhood $U$ biholomorphic to
$\Delta^{n-4}\times V^4$ where $\Delta^{n-4}$ is a complex disc of dimension
$n-4$ and $V^4$ is a neighbourhood of the vertex of the cone over
Veronese emmbedding $\P^1\times\P^2\hookrightarrow\P^5$.
In particular, outside $f^{-1}(Z')$ any non-trivial fiber of $\f$
is $\P^2$ with the normal $\O(-1)^{\oplus 2}\oplus\O^{n-4}$.

\proof The estimate on the dimension of $Z_k$ is obvious. The irreducibility
of $E$ in (i) as well the description of $\f$ outside of
$\f^{-1}(Z_2)$ is in [1], whereas a general blow-up statement is in
[3].  The description of $\f$ in (ii) follows easily from
[11].

\bigskip

For some of the applications of extremal rays which are mentioned in
the introduction it is convenient to know the scheme parametrizing
rational curves from an extremal ray. Below, we sketch the
construction of such an object. For an overview on the theory of
$Hilb(X)$, $Chow(X)$ and $Hom(\P^1,X)$ --- and how they can be used
to parametrize rational curves --- we refer the reader to [13],
Chapter II.1. Rational curves were firstly discussed in this context
by Mori in [15] who proved their existence in extremal rays in
[16].

Let us recall that an {\sl extremal rational curve} in a Mori ray $R$
is a rational curve $C\subset X$ with a normalization $f:\P^1\ra
C\subset X$ such that $[C]\in R$ and its degree
$deg(f^*(-K_X))=-K_X\cdot [C]$ is minimal among the curves from $R$.
In particular, if $R$ is L-negligible then $-K_X\cdot [C]=1$.

Let us take a component ${\cal M}\subset Hom(\P^1,X)$ which contains
the class of $f$. We consider the image of ${\cal M}$ under the
natural map $Hom(\P^1,X)\ra Chow(X)$ which sends $f\in Hom(\P^1,X)$
to the class of $f(\P^1)$. The image of this map (after a
normalisation) is the geometric quotient ${\cal M}/G$ where
$G=Aut(\P^1)$ acts on ${\cal M}$ by composition $G\times{\cal
M}\ni (g,f)\mapsto f\circ g^{-1}\in {\cal M}$.  Now we take the fiber
product ${\cal W}:={\cal M}\times_G \P^1$ which has a natural
projection $q: {\cal W}\ra{\cal M}/G$ and the evaluation $ev: {\cal
W}\ra X$ such that $ev(f,t)=f(t)$.  The quotient ${\cal M}/G$ is a
complete scheme (compact analytic space) and it coincides with the
normalisation of the appropriate component of $Chow(X)$ which
contains $C=f(\P^1)$.  Moreover, at a generic point it is the same as
the appropriate component of $Hilb(X)$. The dimension of ${\cal M}/G$
is bounded from below by Riemann-Roch: $$dim_{\C}{\cal M}/G\geq
-K_C\cdot [C]+n-3.$$
We will say that $C$ has deformations of the
{\sl expected dimension} if the the above inequality becomes an
equality for any component ${\cal M}\subset Hom(\P^1,X)$ which
contains $[f]$.


The line of arguments presented above is due to Mori [15].
We note also the following useful observation

\proclaim Mori Breaking Lemma.
For any two different points $x_1,\
x_2\in X$ the intersection
$q(ev^{-1}(x_1))\cap q(ev^{-1}(x_2))$ is finite.

For L-negligible rays we have the following

\proclaim Lemma 3.5.
Assume that $R$ is an L-negligible Mori ray of a manifold $X$.  If
$C\subset X$ is an extremal rational curve from $R$ then its deformations
are the of expected dimension.

\proof
We use the notation introduced above.  By Breaking Lemma the
dimension of fibers of $(\f_{R})_{|ev({\cal W})}$ is at least by 1
bigger than the dimension of the respective fibers of $ev$.  Indeed,
for any $x\in ev({\cal W})$ the map $ev$ on $q^{-1}(q(ev^{-1}(x)))$
is finite-to-one outside of $ev^{-1}(x)$.
Therefore
$$dim_\C ev({\cal W})-dim_\C \f_R(ev({\cal W}))
\geq dim_\C{\cal W}- dim_\C ev({\cal W})+1.$$
Since $dim_\C{\cal W}=dim_\C {\cal M}/G+1$, in view of 3.1 we have
$$n\geq 2dim_\C ev({\cal W})-dim_\C\f_R(ev({\cal W}))\geq dim_\C{\cal
M}/G+2$$
which is what we want.

\medskip

\noindent{\bf Remark.} Even in the case of a divisorial contraction,
c.f.~3.4, although the exceptional locus of $\f$ is irreducible, the
scheme of extremal rational curves may be reducible. In fact, in 6.9
of [3] we describe a three component Hilbert scheme of extremal
rational curves of a divisorial contraction of a smooth 4-fold with a
reducible fiber which is a degenerated quadric $\P^2\cup\P^2$.  On
the other hand [11] gives an example of a small contraction
of a 4-fold with a disconnected exceptional locus which is a union of
arbitrary (finite) number of disjoint copies of $\P^2$.

We note however that Lemma 3.5 refers to {\sl any} component of the
scheme of extremal rational curves. Moreover, from the above proof of
3.5 it follows that the conclusion of Lemma 3.3 holds also when we
set $S=ev({\cal W})$.  Indeed, in the quoted example from [3] the
special fiber of the divisorial contraction consists of two $\P^2$.

\medskip

Before stating a corollary to the above lemma let us recall that
given a curve $C\subset X$ we say that the class of $C$ remains
holomorphic for small deformations of $X$ if for any smooth family
$\pi: {\cal X}\ra \Delta^1$ of compact complex manifolds over a disc
$\Delta^1$ such that $\pi^{-1}(0)\iso X$ and  any $t$ sufficiently small
there exists a holomorphic curve $C_t\subset\pi^{-1}(t)$ such that
$[C_t]=[C]$ under the natural identification
$H_2(\pi^{-1}(t),\Z)=H_2(X,\Z)$.

\proclaim Corollary 3.6. The class of an extremal rational curve in
an L-negligible Mori ray on a projective manifold $X$
remains holomorphic for small deformations of $X$.

In fact the Corollary is a consequence of 3.5 if we use the following
observation (c.f.~[2]).

\proclaim Lemma 3.7. Let $C\subset X$ be a rational curve in a compact
complex manifold. If the dimension of deformations of $C$ inside $X$
is of the expected dimension then the class of [C] remains
holomorphic for small deformations of $X$.

\proof The dimension of deformations of $C$ inside ${\cal X}$ is bounded
from below by $dim_C{\cal X}-K_{\cal X}\cdot [C]-3=n-K_X\cdot
[C]-2$ so it is bigger than the dimension of deformations of $C$ in
$X$. Therefore the curve $C$ in ${\cal X}$ has to move out from
$\pi^{-1}(0)$ to the neighbouring fibers.


\bigskip

In the conclusion let me make some comments on the geometry of
L-negligible rays in the context of symplectic geometry.  As before,
let us choose an extremal rational curve $C$ and an irreducible
component ${\cal M}\subset Hom(\P^1,X)$ containing the normalisation
of $C$.  For any positive $k$ we consider the product ${\cal
W}_k:={\cal M}\times_G (\P^1)^{\times k}$ where $G=Aut(\P^1)$ acts on
the $k$-th product $(\P^1)^{\times k}$ coordinatewise
$G\times(\P^1)^{\times k}\ni(g,(t_1,\dots t_k))\mapsto (g(t_1),\dots
g(t_k)\in(\P^1)^{\times k}$. On ${\cal W}_k$ we have the natural
projection $q:{\cal W}_k\ra {\cal M}/G$ and the evaluation $ev_k: {\cal
W}_k\ra X^{\times k}$ with $ev_k(f,t_1,\dots t_k)=(f(t_1),\dots
f(t_k))$.  We recall that since $C$ is extremal it follows that
${\cal M}/G$ is compact and $ev_k$ is proper.

We are interested in case $k=2$. In this case however, because of the
breaking lemma, the map $ev_2$ does not have positive dimensional
fibers outside of the diagonal of $X\times X$.  
Thus, if $S_2\subset X^{\times 2}$ is the image of $ev_2$,
then $dim_{\C}S_2=dim_\C{\cal M}-1\geq n-K_X\cdot [C]-2$.
Suppose that $-K_X\cdot [C]=1$ and ${\cal M}$ is of the expected
dimension, then $dim_\C S_2=n$. Let $p_i: X\times X\ra X$ be
the projection on the $i$-factor. Given two cohomology classes
$\alpha_i\in H^{a_i}(X,\Z)$ with $a_1+a_2=2n$ we consider the product
$$(\alpha_1,\alpha_2)_{\cal M}:= (p_1^*(\alpha_1)\cup
p_2^*(\alpha_2))\cdot [S_2].$$

Let $S= p_i(S^2)= ev({\cal M}\times_G\P^1)$ be the locus of curves
from ${\cal M}$ in $X$ and set $s=dim_\C S$.  Moreover for $x\in S$
let $S_x:=ev(q^{-1}(q(ev^{-1}(x))))$. We note that if for some $i=1,\
2$ we have $a_i>2s$ then $\alpha_i\cap p_{i*}(S^2)=0$ and thus
$(\alpha_1,\alpha_2)_{\cal M}=0$.  Thus the product can be non-zero
only if $2s\geq a_i \geq 2n-2s$ for $i=1,\ 2$.  On the other hand if
$a_1=2s$ then we find out that $$(\alpha_1,\alpha_2)_{\cal M}=
(\alpha_1\cdot [S])\cdot(\alpha_2\cdot [S_x])$$ where $x\in S$ is
general.

In particular, if $R$ is a L-negligible extremal ray then the above
discussion applies to any extremal rational curve in $R$. Moreover,
because of the remark following 3.5, for a general $x\in S$ the
normalization of $S_x$ is $\P^{n-s}$.

Similar arguments were used in [18] to compute Gromov-Witten
invariants of some special extremal rational curves. I must admit
however that in the discussion presented above I ignored verification
of the assumptions which are usually asked in the symplectic set-up
to make the computation of Gromov-Witten invariants legible. This can
be done in low dimensional cases by Theorem (3.4).


\beginsection Appendix: remarks and questions.

\par
First let us note
that if $X$ is not projective then the notion of the Lefschetz discriminant
may become void. For example: a Hopf surface.
\medskip
We noted that Hard Lefschetz Theorem implies that $\Delta$ does not meet
the interior ${\cal P}'$ of the nef cone in rational points.
Can $\Delta$ meet ${\cal P}'$ at non-rational points?
\medskip
The bound on the number of L-supported rays obtained by adding degrees
of $\Delta_i$ is probably not the best because subsequent loci $\Delta_i$
are related one to another. What would be the best bound for the number of
L-supported rays?

\medskip
The definitions and the results from Section 1 can be introduced for
complex cohomology ring $H^*(X,\C)$. In particular we can introduce
complex Lefschetz discriminants $\Delta_\C(X)$ and respective
L-supported rays.  It is clear that the ``unbalanced'' components of
the Hodge decomposition of $H^2(X,\C)$, i.e.~$H^{0,2}(X)$ and
$H^{2,0}(X)$ are contained in $\Delta_\C(X)$. Moreover from the
results of Section 2 it follows that if $R$ is an extremal
L-supported ray on $X$ then $R_\C^\perp\supset H^{2,0}(X)\oplus
H^{0,2}(X)$. Indeed, if $\f_R:X\ra Y_R$ is the contraction of $R$
then $H^2(X,\O_X)=\f^*_R H^2(Y,\O_Y)\subset R_\C^\perp$ and obviously
$R_\C^\perp$ is invariant of the complex conjugation of $H^2(X,\C)$

Suppose that $X'$ is another complex structure which defines the
respective Hodge decomposition of $H^2(X,\C)$. Is it true that still
$H^{2,0}(X')\oplus H^{0,2}(X')\subset R_\C^\perp$? What if $X'$ is
obtained by a complex deformation of $X$?  The results of [23] may
suggest that the answer for the second question is positive.

\medskip
The inequality $2dim_\C S-dim_\C\f_R(S)>n$, which appears in
Proposition 3.1, is a sufficient condition for an extremal ray $R$ to
be L-supported.  Let us note however that it destroys the Lefschetz
duality ($\dagger\dagger$) on an even cohomology group of $X$.  I
have been unable to find an example of an extremal ray such that the
inequality $2dim_\C S-dim_\C\f_R(S)\leq n$ is satisfied for any complex
subset $S$ of $X$ and despite of this it is L-supported (e.g.~the
Lefschetz duality fails on odd cohomology).  However it is hard to
expect that the condition from Proposition 3.1 on an extremal ray to
be L-supported is a necessary one.


\beginsection References.

{\smalltype{

\item{[1]} Ando, T., On extremal rays of the higher dimensional
varieties, Invent. Math. {\bf 81} (1985), 347---357.

\item{[2]} Andreatta, M., Peternell, Th.,
On the limits of manifolds with nef canonical bundle,
preprint (1996).

\item{[3]} Andreatta, M., Wi\'sniewski J. A.,
On contractions of smooth varieties,  to appear in J. Alg. Geom.
(199?)

\item{[4]} Andreatta, M., Ballico, E., Wi\'sniewski, J. A., Two theorems on
elementary contractions, Math. Ann. {\bf 297} (1993), 191-198.

\item{[5]} Clemens, H., Koll\'ar, J., Mori, S., {\it Higher dimensional
complex geometry}, Asterisque {\bf 166} (1988).

\item{[6]} Elencwajg, G., The Brauer group in complex geometry, in {\it
Brauer groups in ring theory and algebraic geometry, Antwerp 1981},
Springer Lect.~Notes in Math {917} (1982), 222-230.

\item{[7]} Fujita, T., On polarized manifolds whose adjoint bundles are
not numerically effective, Adv. Stud. Pure Math. {\bf 10} (1987), 167-178.

\item{[8]} Griffiths, Ph., Harris, J., {\it Principles of Algebraic Geometry},
Wiley Interscience Publication, 1978.

\item{[9]} Ionescu, P., Generalized adjunction and applications,
Math.~Proc.~Cambridge Math.~Soc.~{\bf 99} (1988), 457--472.

\item{[10]} Kawamata, Y., Crepant blowing-up of 3-dimensional
canonical singularities and ist application to degeneration of
surfaces, Ann. Math. {\bf 127} (1988), 93--163.

\item{[11]} Kawamata, Y., Small contractions of four dimensional
algebraic manifolds, Math. Ann., {\bf 284} (1989), 595-600.

\item{[12]} Kawamata, Y., Matsuda, K., Matsuki, K.: Introduction to the
Minimal Model Program in {\it Algebraic Geometry, Sendai}, Adv. Studies
in Pure Math. {\bf 10}, Kinokuniya--North-Holland 1987, 283---360.

\item{[13]} Koll\'ar, J., {\it Rational Curves on Algebraic Varieties},
Springer Verlag, Ergebnisse der Math. {\bf 32}, 1995.

\item{[14]} Looijenga, E., Cohomology and intersection homology of algebraic
varieties, in {\it Complex Alg. Geom.} ed. Koll\'ar, AMS-IAS 1997.

\item{[15]} Mori, S., Projective manifolds with ample tangent bundle,
Ann. Math. {\bf 110} (1979), 593--606.

\item{[16]} Mori, S, Threefolds whose canonical bundles are not numerically
effective, Ann. Math., {\bf 116} (1982), 133-176.

\item{[17]} McDuff, D., Salamon, D, {\it J-holomorphic curves and
quantum cohomology}, Univ. Lect. Series {\bf 6} AMS 1994.

\item{[18]} Paoletti, R., On symplectic invariants of algebraic varieties
coming from crepant contractions, preprint (1996)

\item{[19]} Ruan, Y., Symplectic topology and extremal rays,
Geom. Funct. Analysis {\bf 3} (1992) 279--291.

\item{[20]} Wilson, Symplectic deformations of Calabi-Yau threefolds,
preprint (1995)

\item{[21]}  Wi\'sniewski J. A., Length of extremal rays and
generalized adjunction, Math. Zeit. {\bf 200} (1990), 409--427.

\item{[22]} Wi\'sniewski, J. A., On contractions of extremal rays of Fano
manifolds, Jounal f\"ur die reine und angew. Mathematik,
{\bf 417} (1991), 141-157.

\item{[23]} Wi\'sniewski, J. A., On deformation of nef values, Duke Math. J.
{\bf 64} (1991), 325--332.

}}


\bigskip
{\obeylines{
Instytut Matematyki Uniwersytetu Warszawskiego
Banacha 2, 02-097 Warszawa, Poland
e-mail:{\tt jarekw@mimuw.edu.pl}
}}

\end